\documentclass[11pt]{article}

\usepackage{amsmath, epsfig, cite, lineno}
\usepackage{amssymb}
\usepackage{amsfonts, color}
\usepackage{latexsym}

\newtheorem{thm}{Theorem}[section]

\newtheorem{lem}[thm]{Lemma}
\newtheorem{prob}[thm]{\it Problem}
\newcommand{\pf}{\noindent{\it Proof.} }

\def\k{{\bf k}}
\def\z{{\bf z}}
\def\a{{\bf a}}
\def\b{{\bf b}}

\def\0{{\bf 0}}
\def\1{{\bf 1}}
\def\n{{\bf n}}
\def\j{{\bf j}}
\def\r{{\bf r}}
\def\x{{\bf x}}

\numberwithin{equation}{section}

\newcommand{\qed}{{\hfill\rule{4pt}{7pt}}\medskip}

\begin{document}


\begin{center}
{\Large\bf Simple proofs of Jensen's, Chu's, Mohanty-Handa's, and\\[5pt] Graham-Knuth-Patashnik's identities }
\end{center}

\vskip 2mm \centerline{Victor J. W. Guo}

\begin{center}
{\footnotesize Department of Mathematics, East China Normal
University, Shanghai 200062,
 People's Republic of China\\
{\tt jwguo@math.ecnu.edu.cn,\quad
http://math.ecnu.edu.cn/\textasciitilde{jwguo}}}
\end{center}


\vskip 0.7cm
\noindent{\bf Abstract.} Motivated by the recent work of Chu [Electron. J. Combin. 17  (2010), \#N24],
we give simple proofs of Jensen's identity
$$
\sum_{k=0}^{n}{x+kz\choose k}{y-kz\choose n-k}
=\sum_{k=0}^{n}{x+y-k\choose n-k}z^k,
$$
and Chu's and Mohanty-Handa's generalizations of Jensen's identity. We also give a quite simple proof of
an equivalent form of Graham-Knuth-Patashnik's identity
$$
\sum_{k\geq 0}{m+r\choose m-n-k}{n+k\choose n}x^{m-n-k}y^k
=\sum_{k\geq 0}{-r\choose m-n-k}{n+k\choose n}(-x)^{m-n-k}(x+y)^k,
$$
which was rediscovered, respectively, by Sun in 2003 and Munarini in 2005.
Finally we give a multinomial coefficient generalization of this identity
and raise two open problems.


\section{Introduction}
Abel's identity (see, for example, \cite[\S3.1]{Comtet})
\begin{align}
\sum_{k=0}^{n}{n\choose k}x(x+kz)^{k-1}(y-kz)^{n-k}=(x+y)^n \label{eq:abel}
\end{align}
and Rothe's identity (or called Hagen-Rothe's identity, see, for example, \cite[\S5.4]{GKP})
\begin{align}
\sum_{k=0}^{n}\frac{x}{x-kz}{x-kz\choose k}{y+kz\choose n-k}=
{x+y\choose n}, \label{eq:rothe}
\end{align}
are famous in the literature and play an important role in enumerative combinatorics.
Recently, Chu \cite{Chu2010} gave elementary proofs of Abel's identity and Rothe's identity by
using the binomial theorem and the Chu-Vandermonde convolution formula respectively.

Motivated by Chu's work, we shall study Jensen's identity \cite{Jensen}, which is closely related to Rothe's identity,
 and can be stated as follows:
\begin{align}
\sum_{k=0}^{n}{x+kz\choose k}{y-kz\choose n-k}
=\sum_{k=0}^{n}{x+y-k\choose n-k}z^k.
\label{eq:jensen}
\end{align}
Jensen's identity \eqref{eq:jensen} has ever attracted much attention
by different authors. Gould \cite{Gould60} obtained the
following Abel-type analogue:
\begin{align}
\sum_{k=0}^{n}\frac{(x+kz)^k}{k!}\frac{(y-kz)^{n-k}}{(n-k)!}
=\sum_{k=0}^{n}\frac{(x+y)^k}{k!}z^{n-k}.
\label{eq:gould-jensen}
\end{align}
Carlitz \cite{Carlitz} gave two interesting theorem related to \eqref{eq:jensen} and \eqref{eq:gould-jensen} by
mathematical induction.
With the help of generating functions, Gould \cite{Gould62} derived the following variation of Jensen's identity \eqref{eq:jensen}:
\begin{align}
\sum_{k=0}^{n}{x+kz\choose k}{y-kz\choose n-k}=
\sum_{k=0}^{n}k{x+y-k\choose n-k}\frac{x+y-(n-k)z-k}{x+y-k}z^k.
\label{eq:gouldv}
\end{align}
E. G.-Rodeja F. \cite{Rodeja} deduced Gould's identity \eqref{eq:gould-jensen} from \eqref{eq:jensen} by establishing
an identity which includes both. Cohen and Sun \cite{CS} also gave an expression which unifies \eqref{eq:jensen} and \eqref{eq:gould-jensen}.
Chu \cite{Chu86} generalized Jensen's identity \eqref{eq:jensen} to a multi-sum form:
\begin{align}
\sum_{k_1+\cdots+k_{s}=n}\prod_{i=1}^s{x_i+k_i z\choose k_i}
=\sum_{k=0}^{n}{k+s-2\choose k}{x_1+\cdots+x_{s}+nz-k\choose n-k}z^k.\label{eq:chu}
\end{align}
Moreover, the identities \eqref{eq:jensen} and \eqref{eq:chu} were respectively generalized by Mohanty and Handa \cite{MH}
and Chu \cite{Chu89} to the case of multinomial coefficients (to be stated in Section~\ref{sec4}).

The first purpose of this paper is to give simple proofs of Jensen's identity, Chu's identity \eqref{eq:chu},
Mohanty-Handa's identity, and Chu's generalization of Mohanty-Handa's identity.  We shall use the
Chu-Vandermonde convolution formula
\begin{align*}
\sum_{k=0}^n{x\choose k}{y\choose n-k}={x+y\choose n}  
\end{align*}
and a well-known identity
\begin{align}
\sum_{k=0}^{n}(-1)^{n-k}{n\choose k}k^r
=
\begin{cases}
0, &\text{if $0\leq r\leq n-1$,} \\
n!, &\text{if $r=n$.}
\end{cases}  \label{eq:stirling}
\end{align}
Eq. \eqref{eq:stirling} may be easily deduced from the Stirling numbers of the second kind \cite[p.~34, (24a)]{Stanley97}. The first case of
\eqref{eq:stirling} was already utilized by the author \cite{Guo} to give a simple proof of Dixon's identity and by
Chu \cite{Chu2010} in his proofs of Abel's and Rothe's identities.

It is interesting that our proof of Chu's identity \eqref{eq:chu} will also leads to a very short
proof of Graham-Knuth-Patashnik's identity,
which was rediscovered several times in the past few years.
The second purpose of this paper is to give a multinomial coefficient generalization of Graham-Knuth-Patashnik's identity
and raise two open problems.

\section{Proof of Jensen's identity }

By the Chu-Vandermonde convolution formula, we have
\begin{align}
\sum_{k=0}^{n}{x+kz\choose k}{y-kz\choose n-k}
=\sum_{k=0}^{n}{x+kz\choose k}\sum_{i=k}^n{x+y+1\choose n-i}{-x-kz-1\choose i-k} \label{eq:first}
\end{align}
Interchanging the summation order in \eqref{eq:first} and noticing that
$$
{x+kz\choose k}{-x-kz-1\choose i-k}=(-1)^{i-k}{i\choose k}{x+kz+i-k\choose i},
$$
we have
\begin{align}
\sum_{k=0}^{n}{x+kz\choose k}{y-kz\choose n-k}
&=\sum_{i=0}^{n}{x+y+1\choose n-i}\sum_{k=0}^i (-1)^{i-k}{i\choose k}{x+kz+i-k\choose i} \nonumber\\
&=\sum_{i=0}^{n}{x+y+1\choose n-i}(z-1)^i,  \label{eq:new}
\end{align}
where the second equality holds because ${x+kz+i-k\choose i}$ is a polynomial in $k$ of degree $i$
with leading coefficient $(z-1)^i/i!$ and we can apply \eqref{eq:stirling} to simplify.
We now substitute $x\to -x-1$, $y\to -y+n-1$ and $z\to -z+1$ in \eqref{eq:new} and observe that
\begin{align}
{-x\choose k}=(-1)^k{x+k-1\choose k}.\label{eq:negative}
\end{align}
Then we obtain
\begin{align}
\sum_{k=0}^{n}{x+kz\choose k}{y-kz\choose n-k}
=\sum_{i=0}^{n}{x+y-i\choose n-i}z^i,  \label{eq:new2}
\end{align}
as desired.

Combining \eqref{eq:jensen} and \eqref{eq:new}, we get the following identity:
$$
\sum_{k=0}^{n}{x-k\choose n-k}z^k=\sum_{k=0}^{n}{x+1\choose n-k}(z-1)^k,
$$
which is equivalent to the following identity in
Graham et al. \cite[p.~218]{GKP}:
\begin{align*}
\sum_{k\leq m}{m+r\choose k}x^ky^{m-k}=\sum_{k\leq m}{-r\choose k}(-x)^k(x+y)^{m-k}.
\end{align*}

\section{Proofs of Chu's and Graham-Knuth-Patashnik's identities}

By repeatedly using the Chu-Vandermonde convolution formula, we have
\begin{align}
{x_{s}+k_{s} z\choose k_{s}}&={x_{s}+(n-k_1-\cdots-k_{s-1})z\choose n-k_1-\cdots-k_{s-1}}  \nonumber\\
&=\sum_{j=k_1+\cdots+k_{s-1}}^n \sum_{j_1+\cdots+j_{s-1}=j}{x_1+\cdots+x_{s}+nz+s-1\choose n-j}\nonumber\\
&\qquad\ \times \prod_{i=1}^{s-1}{-x_i-k_i z-1\choose j_i-k_i}.
\label{eq:multi-chu}
\end{align}
It follows that 
\begin{align}
\sum_{k_1+\cdots+k_{s}=n}\prod_{i=1}^s{x_i+k_i z\choose k_i}
&=\sum_{k_1+\cdots+k_{s-1}=0}^n \sum_{j=k_1+\cdots+k_{s-1}}^n \sum_{j_1+\cdots+j_{s-1}=j}{x_1+\cdots+x_{s}+nz+s-1\choose n-j}\nonumber\\
&\qquad\ \times \prod_{i=1}^{s-1}{x_i+k_i z\choose k_i}{-x_i-k_i z-1\choose j_i-k_i}.  \label{eq:multisum}
\end{align}
Interchanging the summation order in \eqref{eq:multisum} and observing that
$$
{x_i+k_i z\choose k_i}{-x_i-k_i z-1\choose j_i-k_i}=(-1)^{j_i-k_i}{j_i\choose k_i}{x_i+k_i z+j_i-k_i\choose j_i}
$$
and ${x_i+k_i z+j_i-k_i\choose j_i}$ is a polynomial in $k_i$ of degree $j_i$ with leading coefficient $(z-1)^{j_i}/j_i!$,
by \eqref{eq:stirling} we get
\begin{align}
\sum_{k_1+\cdots+k_{s}=n}\prod_{i=1}^s{x_i+k_i z\choose k_i}
&= \sum_{j=0}^n{x_1+\cdots+x_{s}+nz+s-1\choose n-j} \sum_{j_1+\cdots+j_{s-1}=j}(z-1)^j \nonumber \\
&=\sum_{j=0}^n {j+s-2\choose j}{x_1+\cdots+x_{s}+nz+s-1\choose n-j} (z-1)^j. \label{eq:finalpf}
\end{align}
Substituting $x_i\to -x_i-1$ ($i=1,\ldots,s$) and $z\to -z+1$ in \eqref{eq:finalpf} and using \eqref{eq:negative},
we immediately get Chu's identity \eqref{eq:chu}.

Comparing \eqref{eq:chu} with \eqref{eq:finalpf} and replacing $s$ by $s+2$, we immediately get
\begin{align}
\sum_{k=0}^{n}{k+s\choose k}{x-k\choose n-k}z^k
=\sum_{j=0}^n {k+s\choose k}{x+s+1\choose n-k} (z-1)^k. \label{eq:ks-2}
\end{align}
It is easy to see that the identity \eqref{eq:ks-2} is equivalent to each of the following known identities:
\begin{itemize}
\item Graham-Knuth-Patashnik's identity \cite[p.~218]{GKP}
\begin{align*}
\sum_{k\geq 0}{m+r\choose m-n-k}{n+k\choose n}x^{m-n-k}y^k
=\sum_{k\geq 0}{-r\choose m-n-k}{n+k\choose n}(-x)^{m-n-k}(x+y)^k.
\end{align*}
\item Sun's identity \cite{Sun}
\begin{align}
\sum_{k=0}^m(-1)^{m-k}{m\choose k}{n+k\choose a}(1+x)^{n+k-a}
=\sum_{k=0}^n {n\choose k}{m+k\choose a}x^{m+k-a}.\label{eq:sun}
\end{align}
\item Munarini's identity \cite{Munarini}
\begin{align}
\sum_{k=0}^n(-1)^{n-k}{\beta-\alpha+n\choose n-k}{\beta+k\choose k}(1+x)^k =\sum_{k=0}^n{\alpha\choose n-k}{\beta+k\choose k}x^k.
\label{eq:munarini}
\end{align}
\end{itemize}
Moreover, the following special case
\begin{align}
\sum_{k=0}^n(-1)^{n-k}{n\choose k}{n+k\choose k}(1+x)^k=\sum_{k=0}^n{n\choose k}{n+k\choose k}x^k
\label{eq:simons}
\end{align}
was reproved by Simons \cite{Simons}, Hirschhorn \cite{Hirschhorn}, Chapman \cite{Chapman}, Prodinger
\cite{Prodinger}, Wang and Sun \cite{WS}.

\section{Mohanty-Handa's identity and Chu's generalization\label{sec4}}
Let $m$ be a fixed positive integer.
For $\a=(a_1,\ldots,a_m)\in\mathbb{N}^m$ and
$\b=(b_1,\ldots,b_m)\in\mathbb{C}^m$, set $|\a|=a_1+\cdots+a_m$,
$\a!=a_1!\cdots a_m!$,
$\a+\b=(a_1+b_1,\ldots,a_m+b_m)$, $\a\cdot\b=a_1 b_1+\cdots+a_m b_m$,
and $\b^\a=b_1^{a_1}\cdots b_m^{a_m}$.
For any variable $x$ and $\n=(n_1,\ldots,n_m)\in\mathbb{Z}^m$,
the {\it multinomial coefficient} ${x\choose \n}$ is defined by
\begin{align*}
{x\choose \n}=
\begin{cases}
x(x-1)\cdots(x-|\n|+1)/\n!,
&\text{if $\n\in\mathbb{N}^m$,}\\
0, &\text{otherwise.}
\end{cases}
\end{align*}
Moreover, we let $\0=(0,\ldots,0)$ and $\1=(1,\ldots,1)$.

In 1969, Mohanty and Handa \cite{MH} established the following multinomial coefficient generalization of Jensen's identity
\begin{align}
\sum_{\k=\0}^{\n}{x+\k\cdot\z\choose \k}{y-\k\cdot\z\choose \n-\k}
=\sum_{\k=\0}^{\n}{x+y-|\k|\choose \n-\k}{|\k|\choose \k}\z^{\k}.
\label{eq:mh}
\end{align}
Twenty years later,  Mohanty-Handa's identity was generalized by Chu \cite{Chu89} as follows:
\begin{align}
\sum_{\k_1+\cdots+\k_{s}=\n}\prod_{i=1}^s{x_i+\k_i\cdot \z\choose \k_i}
=\sum_{\k=\0}^{\n}{|\k|+s-2\choose \k}{x_1+\cdots+x_{s}+\n\cdot\z-|\k|\choose \n-\k}\z^\k, \label{eq:chu89}
\end{align}
which is also a generalization of \eqref{eq:chu}.

\medskip
\noindent{\it Remark.}
Note that the corresponding multinomial coefficient generalization of Rothe's identity was already obtained
by Raney \cite{Raney} (for a special case) and Mohanty \cite{Mohanty66}.
The reader is referred to Strehl \cite{ Strehl} for a historical note on Raney-Mohanty's identity.
\medskip

In what follows, we will give an elementary proof of Chu's identity \eqref{eq:chu89} similar to that of \eqref{eq:chu}.
First note that the Chu-Vandermonde convolution formula has the following trivial generalization
\begin{align}
\sum_{\k=\0}^\n{x\choose \k}{y\choose \n-\k}={x+y\choose \n},  \label{eq:cvmulti}
\end{align}
as mentioned by Zeng \cite{Zeng}, while \eqref{eq:stirling} can be easily generalized as
\begin{align}
\sum_{\k=\0}^{\n}(-1)^{|\n|-|\k|}{\n\choose \k}\k^\r
=
\begin{cases}
0, &\text{if $r_i<n_i$ for some $1\leq i\leq m$.} \\
\n!, &\text{if $\r=\n$,}
\end{cases}  \label{eq:multi-stirling}
\end{align}
where
$$
{\n\choose \k}:=\prod_{i=1}^m{n_i\choose k_i}.
$$

\begin{lem}For $\n\in\mathbb{N}^m$ and $s\geq 1$, there holds
\begin{align}
\sum_{\k_1+\cdots+\k_{s}=\n}\prod_{i=1}^s{|\k_i|\choose \k_i}
={|\n|+s-1\choose \n}. \label{eq:lem}
\end{align}
\end{lem}

\pf For nonnegative integers $a_1,\ldots, a_s$ such that $a_1+\cdots+a_s=|\n|$, by the Chu-Vandermonde convolution formula \eqref{eq:multi-stirling},
the following identity holds
\begin{align}
\sum_{\k_1+\cdots+\k_{s}=\n}\prod_{i=1}^s{a_i\choose \k_i}
={|\n|\choose \n}. \label{eq:multicv2}
\end{align}
Moreover, in this case, for $\k_1+\cdots+\k_{s}=\n$, we have
$$
\prod_{i=1}^s{a_i\choose \k_i}\neq 0\quad \text{if and only if}\quad |\k_i|=a_i\ (i=1,\ldots,s).
$$
Thus, the identity \eqref{eq:multicv2} may be rewritten as
\begin{align*}
\sum_{\substack{\k_1+\cdots+\k_{s}=\n\\|\k_1|=a_1,\ldots,|\k_s|=a_s}}
\prod_{i=1}^s{a_i\choose \k_i}
={|\n|\choose \n}. 
\end{align*}
It follows that
\begin{align*}
\sum_{\k_1+\cdots+\k_{s}=\n}\prod_{i=1}^s{|\k_i|\choose \k_i}
&=\sum_{a_1+\cdots+a_s=|\n|}\sum_{\substack{\k_1+\cdots+\k_{s}=\n\\|\k_1|=a_1,\ldots,|\k_s|=a_s}}
\prod_{i=1}^s{a_i\choose \k_i}  \\
&=\sum_{a_1+\cdots+a_s=|\n|}{|\n|\choose \n} \\
&={|\n|+s-1\choose |\n|}{|\n|\choose \n},
\end{align*}
as desired.
\qed

By repeatedly using the convolution formula \eqref{eq:cvmulti}, we may rewrite
the left-hand side of \eqref{eq:chu89} as
\begin{align}
&\sum_{\k_1+\cdots+\k_{s-1}=\0}^\n \sum_{\j=\k_1+\cdots+\k_{s-1}}^\n
\sum_{\j_1+\cdots+\j_{s-1}=\j}{x_1+\cdots+x_{s}+\n\cdot\z+m-1\choose \n-\j}\nonumber\\
&\qquad\ \times \prod_{i=1}^{s-1}{x_i+\k_i\cdot\z\choose \k_i}{-x_i-\k_i\cdot \z-1\choose \j_i-\k_i}.  \label{eq:multi-chufi}
\end{align}
Interchanging the summation order in \eqref{eq:multi-chufi}, observing that
$$
{x_i+\k_i\cdot \z\choose \k_i}{-x_i-\k_i\cdot \z-1\choose \j_i-\k_i}
=(-1)^{|\j_i|-|\k_i|}{\j_i\choose \k_i}{x_i+\k_i\cdot \z+|\j_i|-|\k_i|\choose \j_i}
$$
and
$$
{x_i+\k_i\cdot \z+|\j_i|-|\k_i|\choose \j_i}
$$
is a polynomial in $k_{i\,1},\ldots,k_{i\,m}$ with the coefficient of $\k_i^{\j_i}$ being ${|\j_i|\choose \j_i}(\z-\1)^{\j_i}/\j_i!$.
Applying \eqref{eq:multi-stirling}, we get
\begin{align}
&\hskip -3mm \sum_{\k_1+\cdots+\k_{s}=\n}\prod_{i=1}^s{x_i+\k_i\cdot \z\choose \k_i}  \nonumber\\
&= \sum_{\j=\0}^\n{x_1+\cdots+x_{s}+\n\cdot\z+s-1\choose \n-\j} (\z-\1)^\j
\sum_{\j_1+\cdots+\j_{s-1}=\j}\prod_{i=1}^m {|\j_i|\choose \j_i}\nonumber \\
&=\sum_{\j=\0}^\n {|\j|+s-2\choose \j}{x_1+\cdots+x_{s}+\n\cdot\z+s-1\choose \n-\j} (\z-\1)^\j, \label{eq:multi-finalpf}
\end{align}
where the second equality follows from \eqref{eq:lem}.
Substituting $x_i\to -x_i-1$ ($i=1,\ldots,s$) and $\z\to -\z+\1$ in \eqref{eq:multi-finalpf} and observing that
${-x\choose \k}=(-1)^{|\k|}{x+|\k|-1\choose \k}$, we immediately get \eqref{eq:chu89}.

Comparing \eqref{eq:chu89} with \eqref{eq:multi-finalpf} and replacing $s$ by $s+2$, we obtain the following result.
\begin{thm}For $\n\in\mathbb{N}^m$ and $\z\in\mathbb{C}^m$, there holds
\begin{align}
\sum_{\k=\0}^{\n}{|\k|+s\choose \k}{x-|\k|\choose \n-\k}\z^\k
=\sum_{\k=\0}^\n {|\k|+s\choose \k}{x+s+1\choose \n-\k} (\z-\1)^\k. \label{eq:newmulti}
\end{align}
\end{thm}

It is easy to see that \eqref{eq:newmulti} is a multinomial coefficient generalization of \eqref{eq:ks-2}.
Substituting $s\to\beta$, $x\to\alpha-\beta-1$ and $\z\to\1+\x$ in \eqref{eq:newmulti}, we get
\begin{align}
\sum_{\k=\0}^\n(-1)^{|\n|-|\k|}{\beta-\alpha+|\n|\choose \n-\k}{\beta+|\k|\choose \k}(\1+\x)^\k
 =\sum_{\k=\0}^\n{\alpha\choose \n-\k}{\beta+|\k|\choose \k}\x^\k,
\label{eq:multi-munarini}
\end{align}
which is a generalization of Munarini's identity \eqref{eq:munarini}. If $\alpha=\beta=|\n|$, then \eqref{eq:multi-munarini}
reduces to
\begin{align*}
\sum_{\k=\0}^\n(-1)^{|\n|-|\k|}{|\n|\choose \n-\k}{|\n|+|\k|\choose \k}(\1+\x)^\k
 =\sum_{\k=\0}^\n{|\n|\choose \n-\k}{|\n|+|\k|\choose \k}\x^\k,
\end{align*}
which is generalization of Simons' identity \eqref{eq:simons}.
 Note that Shattuck \cite{Shattuck} and  Chen and Pang \cite{CP} have given
different combinatorial proofs of \eqref{eq:munarini}. It is natural to ask

\begin{prob}{\rm
Find a combinatorial interpretation of \eqref{eq:multi-munarini}. }
\end{prob}

\section{Concluding remarks}
We know that binomial coefficient identities usually have
nice $q$-analogues. However, there are only curious (not natural)
$q$-analogues of Abel's and Rothe's identities (see \cite{Schlosser}
and references therein) up to now. There seems to have no
$q$-analogues of Jensen's identity in the literature.

It is interesting that Hou and Zeng \cite{HZ} gave a $q$-analogue of Sun's identity \eqref{eq:sun}:
\begin{align}
\sum_{k=0}^m(-1)^{m-k}{m\brack k}{n+k\brack a}(-xq^a;q)_{n+k-a}q^{{k+1\choose 2}-mk+{a\choose 2}} =\sum_{k=0}^n
{n\brack k}{m+k\brack a}x^{m+k-a}q^{mn+{k\choose 2}},  \label{eq:hz}
\end{align}
where
$(a;q)_n=(1-a)(1-aq)\cdots(1-aq^{n-1})$ and
$$
{\alpha\brack k} =
\begin{cases}
\displaystyle\frac{(q^{\alpha-k+1};q)_k}{(q;q)_k}, &\text{if $k\geq 0$},\\[10pt]
0,&\text{if $k<0$.}
\end{cases}
$$
Clearly, \eqref{eq:hz} may be written as a $q$-analogue of Munarini's identity \eqref{eq:munarini}:
\begin{align}
&\hskip -3mm
\sum_{k=0}^n(-1)^{n-k}{\beta-\alpha+n\brack n-k}{\beta+k\brack k}q^{{n-k\choose 2}-{n\choose 2}}(-x;q)_k \nonumber\\
&=\sum_{k=0}^n{\alpha\brack n-k}{\beta+k\brack k}q^{{n-k+1\choose 2}+(\beta-\alpha)(n-k)} x^k,   \label{eq:hznew}
\end{align}
as mentioned by Guo and Zeng \cite{GZ}. We end this paper with the following problem.

\begin{prob}{\rm
Is there a $q$-analogue of \eqref{eq:multi-munarini}? Or equivalently, is there a multi-sum generalization of \eqref{eq:hznew}?}
\end{prob}

\vskip 2mm
\noindent{\bf Acknowledgments.} This work was partially
supported by the Fundamental Research Funds for the Central Universities,  Shanghai Rising-Star Program (\#09QA1401700),
Shanghai Leading Academic Discipline Project (\#B407), and the National Science Foundation of China (\#10801054).


\begin{thebibliography}{99}
\small \setlength{\itemsep}{-.8mm}

\bibitem{Carlitz}L. Carlitz, Some formulas of Jensen and Gould,
Duke Math. J. 27 (1960), 319--321.

\bibitem{Chapman}R. Chapman, A curious identity revisited, Math. Gazette 87 (2003), 139--141.

\bibitem{CP}W. Y. C. Chen and S. X. M. Pang, On the combinatorics of the Pfaff identity,
 Discrete Math. 309 (2009), 2190--2196.

\bibitem{Chu86}W. Chu, On an extension of a partition identity and its Abel analog,
J. Math. Rese. Exposition 6 (4) (1986), 37--39.

\bibitem{Chu89}W. Chu, Jensen's theorem on multinomial coefficients and its
Abel-analog, Appl. Math. J. Chinese Univ. {4} (1989), 172--178 (in Chinese).

\bibitem{Chu2010}W. Chu, Elementary proofs for convolution identities of Abel and Hagen-Rothe,
Electron. J. Combin. 17  (2010), \#N24.

\bibitem{CS}M. E. Cohen and H. S. Sun, A note on the Jensen-Gould convolutions, Canad. Math. Bull. 23 (1980), 359--361.

\bibitem{Comtet}L. Comtet, Advanced Combinatorics, D. Reidel Publishing Company, Dordrecht-Holland, 1974.

\bibitem{GKP} R. L. Graham, D. E. Knuth and O. Patashnik,
Concrete Mathematics, Addion-Wesley Pubilshing Co., 1989.

\bibitem{Rodeja}E. G.-Rodeja F., On identities of Jensen, Gould and Carlitz,
in: Proc. Fifth Annual Reunion of Spanish Mathematicians (Valencia, 1964),
Publ. Inst. ``Jorge Juan" Mat., Madrid, 1967, pp. 11--14.

\bibitem{Gould60}H. W. Gould, Generalization of a theorem of Jensen concerning convolutions,
Duke Math. J. 27 (1960) 71--76.

\bibitem{Gould62}H. W. Gould, Involving sums of binomial coefficients and a formula of Jensen,
 Amer. Math. Monthly, 69 (5) (1962), 400--402.

\bibitem{Guo}V. J. W. Guo, A simple proof of Dixon's identity, Discrete Math. 268 (2003), 309--310.

\bibitem{GZ}V. J. W. Guo and J. Zeng, Combinatorial proof of a curious $q$-binomial coefficient identity,
Electron. J. Combin. 17 (2010), \#N13.

\bibitem{HZ}S. J. X. Hou and J. Zeng, A $q$-analog of dual sequences with applications,
 European J. Combin. 28 (2007), 214--227.

\bibitem{Hirschhorn}M. Hirschhorn, Comment on a curious identity, Math. Gazette 87 (2003), 528--530.

\bibitem{Jensen}J. L. W. V. Jensen, Sur une identit\'e d'Abel et sur d'autres formules analogues,
Acta Math. 26 (1902), 307--318.

\bibitem{Mohanty66}S. G. Mohanty,
Some convolutions with multinomial coefficients and related probability distributions,
SIAM Rev. {8} (1966), 501--509.

\bibitem{Mohanty}S. G. Mohanty, Lattice Path Counting and Applications,
Academic Press, New York, 1979.

\bibitem{MH}S.G. Mohanty and B.R. Handa,
Extensions of Vandermonde type convolutions with several summations and their applications,
I. Canad. Math. Bull. 12 (1969), 45--62.

\bibitem{Munarini}E. Munarini, Generalization of a binomial identity of Simons, Integers 5 (2005), \#A15.

\bibitem{Prodinger}H. Prodinger, A curious identity proved by Cauchy's integral formula, Math. Gazette 89
(2005), 266--267.

\bibitem{Raney}G. N. Raney, Functional composition patterns and
power series reversion, Trans. Amer. Math. Soc. {94} (1960), 441--451.

\bibitem{Rothe}H. A. Rothe, Formulae de serierum reversione demonstratio
universalis signis localibus combinatorio-analyticorum vicariis exhibita,
Leipzig, 1793.

\bibitem{Schlosser}M. Schlosser, Abel-Rothe type generalizations of
Jacobi's triple product identity, in: Theory and Applications of Special
Functions, Dev. Math., 13, Springer, New York, 2005, pp.~383--400.


\bibitem{Shattuck}M. Shattuck, Combinatorial proofs of some Simons-type binomial coefficient identities,
Integers 7 (2007), \#A27.

\bibitem{Simons}S. Simons, A curious identity, Math. Gazette 85 (2001), 296--298.

\bibitem{Stanley97}R. P. Stanley, Enumerative Combinatorics, Vol. 1,
Cambridge Studies in Advanced Mathematics, 49, Cambridge University Press, Cambridge, 1997.

\bibitem{Strehl}V. Strehl, Identities of Rothe-Abel-Schl\"afli-Hurwitz-type,
Discrete Math. {99}  (1992),  321--340.

\bibitem{Sun}Z.-W. Sun, Combinatorial identities in dual sequences, European J. Combin. 24 (2003), 709--718.

\bibitem{WS}X. Wang and Y. Sun, A new proof of a curious identity, Math. Gazette 91 (2007), 105--106.

\bibitem{Zeng}J. Zeng, Multinomial convolution polynomials,
Discrete Math. {160} (1996), 219--228.

\end{thebibliography}
\end{document}